\documentclass[12pt]{article}      
\usepackage{amssymb}
\newtheorem{Pro}{Proposition}

\begin{document} \sloppy 

\title{Rational cohomology of an algebra need not be detected by Frobenius kernels}

\author{Wilberd van der Kallen}\newcommand{\gr}{{\mathrm{gr}}}
\date{}
\maketitle

\begin{abstract}
We record some cohomological computations in an example of Friedrich Knop.
The example is a polynomial algebra  in 
characteristic two with an unusual $\mathrm {SL}_2$ action.
\end{abstract}

\section{The example}
In \cite[section 5]{Knop} Friedrich Knop gave an example of a transitive action of the algebraic group $\mathrm {SL}_2$ 
on the afine plane $\mathbb A^2$ in characteristic two.
At the 2002 CRM Workshop on Invariant Theory, cf.~\cite{vdK coh}, he also explained that the example has noteworthy 
properties in connection with Grosshans grading. In this note we look at the rational cohomology and see that
it is equally instructive.
Thus consider the algebraic group $G=\mathrm {SL}_2$ defined over a field $k$ of characteristic two.
The diagonal subgroup is $T$, the unipotent upper triangular subgroup is $U$.
Let $N$ be the normalizer of $T$ in $G$. The example is then $G/N\cong \mathbb A^2$.
As Knop observed, $k[G/N]\hookrightarrow k[G/T]$ is separable,
$k[G/T]$ has good filtration and $k[G/N]^U=k[G/T]^U$.
Recall from \cite[2.3]{vdK coh} that the Grosshans graded $\gr\, k[G/T]$, known as the `hull' of $\gr\, k[G/N]$, 
is a purely inseparable extension of
$\gr\, k[G/N]$. So in some sense the inseparability is a property of the Grosshans
filtration, not of the ring extension. We will give cohomology
computations that amplify these observations.

We write the general matrix in $G$ as  $\pmatrix{a&b\cr c&d}$, so that $$k[G]=k[a,b,c,d]/(ad-bc-1).$$
Recall that $k[U\backslash G]=k[c,d]$, a polynomial ring in two variables. Here by $k[U\backslash G]$ we
mean the ring of rational functions on $G$ that are invariant under left translation
by elements of $U$. (One may check $k[U\backslash G]=k[c,d]$
using the multiplicities in a good filtration of $k[G]$ as $G\times G$ module.)
Thus $k[G/T]^U=k[U\backslash G]^T=k[cd]=k[U\backslash G]^N=k[G/N]^U$, a polynomial ring in one variable,
written $cd$ to indicate its image in $k[G]$.
Note that $ad\in  k[G/T]$ satifies $ad(ad-1)=abcd\in k[G/N]$, while no power $(ad)^{2^r}$ is in $k[G/N]$,
because the involution $\sigma=\pmatrix{0&1\cr 1&0}\in N$ interchanges $ad$ with $bc$, but $(ad)^{2^r}-(bc)^{2^r}=1$.
It is easy to see that $k[G/T]$ is generated by $ad$, $ab$, $cd$.
Using the involution $\sigma$ again, one sees that $k[G/T]$ is a free $k[G/N]$-module with basis $1$, $ad$.
Thus $k[G/N]\hookrightarrow k[G/T]$ is separable, $x^2-x-(ab)(cd)$ being the minimal polynomial of
$ad$.  One also finds that $k[G/N]$ is a polynomial ring in two variables, called
$ab$ and $cd$. One may find a basis of $k[G/T]$ and make
its (good) Grosshans filtration explicit as a check of the above.

The Grosshans filtration of  $k[G/T]$ starts with the span of $1$, then the span of $1$, $ab$, $ad$, $cd$.
Intersecting with $k[G/N]$ one gets the span of $1$ and the span of $1$, $ab$, $cd$ respectively.
So the class of $ad$ is an invariant in $k[G/T]/k[G/N]$ and this defines an extension
of $k$ by $k[G/N]$
that does not split because $k[G/N]^G=k[G/T]^G=k$. 
But $ad-(ad)^{2^r}\in k[G/N]$ for all $r\geq1$, so for every such $r$ the extension splits
as an extension of modules for the Frobenius kernel $G_r$, with $(ad)^{2^r}\in k[G/T]^{G_r}$ as lift of
the class of $ad$.
We have shown:
\begin{Pro}
$H^1(G,k[G/N])$ is not the inverse limit of the $H^1(G_r,k[G/N])$.
\end{Pro}
This is in contrast with what one knows for finite
dimensional representations \cite[II 4.12]{J}.

One may go a little further and compute both $H^*(G,k[G/N])$ and $H^*(G_r,k[G/N])$.
As $k[G/T]/k[G/N]$ is a free $k[G/N]$ module on one generator, we get
an the extension $${\cal E}:\quad0\to k[G/N]\to k[G/T]\to k[G/N]\to0.$$
Further  $k[G/T]$ is a direct summand of the injective module $k[G]$,
so it is easy to compute that $H^i(G,  k[G/N])=k$ for all $i\geq0$.
More specifically, if $f$ denotes the inclusion $k\hookrightarrow  k[G/N]$, then $H^i(G,  k[G/N])=k$
is spanned for $i\geq1$ by the Yoneda product of $f$ and $i$ copies of $\cal E$.
So, by \cite[3.2]{Benson I}, $H^*(G,k[G/N])$ is a polynomial ring in one variable of degree one.
 
Let $r\geq1$.
Now $k[G]$ is also $G_r$ injective, by \cite[I 4.12, 5.13]{J}. But $\cal E$ splits over $G_r$, so one
gets $H^i(G_r,k[G/N])=0$ for $i>0$. So the $G_r$ detect none of the higher cohomology of $k[G/N]$.
One may also check that $H^0(G_r,k[G/N])^{(-r)}\cong k[G/N]$. Of course $H^0(G,k[G/N])=k$
is the inverse limit of the $H^0(G_r,k[G/N])$.

We leave it to the reader to compute the nontrivial $H^*(G_1,\gr\; k[G/N])$. (Use \cite[3.10]{Andersen-Jantzen}.)

\end{document}